\documentclass[10pt]{article}
\usepackage{amsmath}
\usepackage{latexsym}
\usepackage{amssymb}
\usepackage{amsfonts}
\usepackage{amsthm}
\title{ARNOLD'S  HYDRODYNAMICS  REVISITED}
\author{J.F. Pommaret \\ CERMICS, Ecole Nationale des Ponts et
  Chauss\'ees,\\6/8 Av. Blaise Pascal, 77455 Marne-la-Vall\'ee Cedex 02,
  France \\
e-mail: pommaret@cermics.enpc.fr \\ (http://cermics.enpc.fr/$\sim$pommaret/home.html)}
\date{  }
\begin{document}
\maketitle

\noindent
{\bf ABSTRACT}:  The purpose of this short paper is to revisit the infinite Lie group theoretical framework of hydrodynamics developped by V. Arnold in 1966. First of all, we extend this approach from the Lie pseudogroup of volume preserving transformations to an arbitrary Lie pseudogroup. Then we prove that, contrary to what could be believed from the work of Arnold which is of a {\it purely analytical} nature, the same results can be obtained from a {\it purely formal} point of view. Finally, we provide the analogue for {\it both} the so-called "{\it body}" and "{\it space}" dynamical equations. We conclude by showing that even this new approach can be superseded by dynamics on Lie groupoids, along ideas pioneered by the brothers E. and F. Cosserat or H. Weyl, on the condition to change the underlying philosophy.\\


\noindent
{\bf INTRODUCTION}:\\

In a celebrated paper published in 1966 [1], V. Arnold applies the differential geometry of infinite dimensional Lie groups to the hydrodynamics of perfect uncompressible fluids. His technique amounts to use for the Lie pseudogroup of volume preserving diffeomorphisms of a bounded Riemannian domain, the analytical analogue of the formal methods already developped by H. Poincar\' e [13], N.G. Chetaev [5,6], E. and F. Cosserat [7,8,9,11], G. Birkhoff [3] and himself [1,2] for the dynamics of a rigid body and, more generally, for the dynamics that can be achieved on any Lie group. It must be noticed that the last reference does not quote the previous ones and we advise the reader to compare [3, p 205 to 216, in particular the formula on p 215/216] with [2,Appendice 2, in particular the formula of Th 2.1 on p 326].\\
The basic tool, which is crucially used, is the so-called "{\it orthogonal decomposition theorem}" of H. Weyl roughly saying that any vector field on the previous domain can be decomposed uniquely into the sum of a divergence-free vector field on this domain, tangent to the boundary, and the gradient of a univalent function.\\
Since the time we read this paper in the seventies, we have always been convinced that this analytical approach could be replaced by a purey formal approach following, for Lie pseudogroups, the one already existing for Lie groups, the latter one being the founding stone of {\it gauge theory} [21]. It is only now that we have been able to succeed while revisiting once more Arnold's paper in the light of new results recently obtained for the {\it partial differential} (PD) {\it optimal control theory} and the corresponding multivariable variational problem with differential constraints [17].\\
We first present a new interpretation of the dynamics on Lie groups in the form of a differential sequence called {\it gauge sequence} and exhibit the corresponding linearized sequence.\\ 
We then extend the preceding results to Lie pseudogroups and obtain the corresponding differential sequence both with the corresponding variational calculus with constraints. All these results are new and the particular case of a one-parameter (time) gauging of the Lie pseudogroup of volume preserving transformations provides the equations of Arnold plus the corresponding "space" equations that he did not obtain. The corresponding second section has been deliberately written in a self-contained way though it relies on tricky calculations involving the implicit function theory.\\
Finally, revisiting the proof of one key theorem within the previous specific case of hydrodynamics, we provide doubts about the usefulness of this approach and hints for using {\it Lie groupoids} and {\it Lie algebroids} instead of {\it Lie pseudogroups} and {\it Lie algebras of vector fields} (For a nice introduction to Lie groupoids, see [12, appendix]). \\
In the case of Lie group actions, studying only the linear framework for simplicity while introducing the Spencer operator, we exhibit the link existing between the linear gauge sequence and the corresponding {\it Spencer sequence} for Lie algebroids [14,15,16,18]. But this is just the generalized Cosserat theory that we already published [14].\\
As global actions on manifolds may not exist, all the results in this paper are local ones though global notations are used for simplicity in order to avoid using explicit open sets and charts.\\

\noindent
{\bf DYNAMICS ON LIE GROUPS}:\\

This section, which is a summary of results already obtained in [14,15], is provided for fixing the notations and the techniques leading to the gauge sequence and the corresponding variational calculus, both with the respective linearized versions.\\
Let $X$ be a manifold of dimension $n$ with local coordinates $x=(x^1,...,x^n)$ and latin indices $i,j=1,...,n$. A point on $X$ will play the part of $n$ {\it parameters}. We denote by $T=T(X)$ the tangent bundle to $X$ and by $T^*=T^*(X)$ the cotangent bundle to $X$ while ${\wedge}^rT^*$ is the bundle of $r$-forms on $X$. Let now $G$ be a Lie group of dimension $p$ with identity $e$, local coordinates $a=(a^1,...,a^p)$ and greek indices $\rho,\sigma,\tau$. We denote by ${\cal{G}}=T_{e}(G)$ the corresponding Lie algebra with vectors denoted by greek letters $\lambda,\mu,\nu$. As usual, we shall identify a map $a:X\rightarrow G$ called {\it gauging} of $G$ over $X$, with its graph $X\rightarrow X\times G$ which is a section of a trivial principal bundle, and, similarly, use the same notation for a bundle and its sheaf of (local) sections as the background will always tell the right choice. In particular, when differential operators are involved, the sectional point of view must automatically be used. Such a convention allows to greatly simplify the notations at the expense of a slight abuse of language.\\
If we have a map $a:X\rightarrow G:x\rightarrow a(x)$, we obtain the tangent map $T(a):T(X)\rightarrow T_a(G):dx\rightarrow da=\frac{\partial a}{\partial x}dx$ and we can pull back the image to $T_e(G)=\cal{G}$ by acting with the inverse $a^{-1}$ of $a$, either on the {\it left} to get $A=a^{-1}da\in T^*\otimes \cal{G}$ or on the {\it right} to get $B=daa^{-1}\in T^*\otimes \cal{G}$. Differentiating $a^{-1}a=e$, we get $a^{-1}da=-ad(a^{-1})$ and thus $daa^{-1}=-ad(a^{-1})$, that is one obtains $B$ from $A$ by changing $a$ to $a^{-1}$ {\it and} changing the sign too. Also we obtain symbolically $A=a^{-1}Ba=Ad(a)B$ by introducing the {\it adjoint map} $\mu=Ad(a)\lambda$ obtained by carrying $\lambda\in \cal{G}$ from $e$ to $a$ on the left and coming back to $\mu\in \cal{G} $ from $a$ to $e$ on the right. For more technical details on the adjoint map, we refer the reader to [15,p 180].\\
Another way to present the previous construction is to introduce the Maurer-Cartan left invariant 1-form $\omega=({\omega}^{\tau}_{\sigma}(a)da^{\sigma})$ on $G$ with value in $\cal(G)$ and pull it back on $X$ by $T(a)$ in order to get $A=A^{\tau}_idx^i={\omega}^{\tau}_{\sigma}(a(x)){\partial}_ia^{\sigma}(x)dx^i$. Finally, using the Maurer-Cartan structure equations for $\omega$ on $G$, namely $d{\omega}^{\tau}-c^{\tau}_{\rho\sigma}{\omega}^{\rho}\wedge{\omega}^{\sigma}=0$, and pulling them back on $X$ similarly, we get a well defined operator $T^*\otimes {\cal{G}} \rightarrow {\wedge}^2T^*\otimes {\cal{G}}:A\rightarrow dA-[A,A]=F$ with local coordinates $F^{\tau}_{ij}={\partial}_iA^{\tau}_j-{\partial}_jA^{\tau}_i-c^{\tau}_{\rho\sigma}A^{\rho}_iA^{\sigma}_j$ where the $c$ are (care to the sign) the structure constants on $\cal{G}$ with Lie algebra bracket $([\lambda,\mu])^{\tau}=c^{\tau}_{\rho\sigma}{\lambda}^{\rho}{\mu}^{\sigma}$. We may collect these results in the following theorem:\\

\noindent
{\bf THEOREM} 1: In the previous framework, there exists the so-called {\it gauge sequence}:\\
\[  \begin{array}{ccccc}
X\times G & \rightarrow  &  T^*\otimes \cal{G}  &  \rightarrow &{\wedge}^2T^*\otimes \cal{G}   \\
   a   &  \rightarrow  &  a^{-1}da=A  &    &      \\
        &                        &   A    &  \rightarrow &   dA-[A,A]=F 
        \end{array}        \]
        
  \noindent
 {\bf REMARK 1}: The previous results can be extended to connections on principal bundles and their curvature [15] but it is important to notice that, in both cases, the group is not acting on the base space.\\
  
 \noindent
 {\bf REMARK 2}: When $n=1$, no differential sequence is existing and this is the situation considered by Arnold for the only "time" parameter $t$. In the case of a rigid body moving in ${\mathbb{R}}^3$, changing slightly the notations with $x_0$ the initial position and $x$ the final position at time $t$, the movement of the rigid body is $x=a(t)x_0+b(t)$ and the {\it projection of the speed in the body} is $a^{-1}\dot{x}=a^{-1}\dot{a}x_0+a^{-1}\dot{b}$ with standard notation for time derivative. This result brings out at once the form $A$ appearing for the Lie group of rigid motions where $G=(a,b)$ with parameters $a$ for rotations and $b$ for translations. Nevertheless, when $n=2$, with parameters the curvilinear abcissa $s$ and the time $t$, the above framework is the one adapted to the Kirchoff theory of thin elastic beams along the work of E. and F. Cosserat for the group of rigid motions [9]. Also, when $n=3$, the two previous operators are exactly described by the brothers Cosserat in the nice reference [11] which has never been quoted elsewhere. An important but tricky question raised by mechanicians was thus to understand why the compatibility conditions of Cosserat theory were first order PD eqations while they were known to be second order PD equations in classical elasticity theory [11]. \\
 
 \noindent
 {\bf REMARK 3}: The second operator has been introduced by E. Cartan for introducing the so-called {\it curvature} (G is the rotation group) and {\it torsion} (G is the translation group), {\it but with no reference to the first operator} [4]. It must be noticed that the "{\it abstract}" {\it group} $G$ {\it is not acting on} $X$. Therefore, if one wants to relate the above framework to electromagnetism (EM), {\it the only possibility} is to consider $X$ as space-time and to call $A$ the 4-EM potential, $F$ the EM field (made by $\vec{E}$ and $\vec{B}$ combined together), with the {\it necessary condition} to have ${dim(\cal{G})}=1$. This was the birth of gauge theory with $G=U(1)$, the unit circle in the complex plane and $\cal{G}$ the parallel to the complex axis at he point $(1,0)$.\\
 
 It just remains to introduce the previous results into a variational framework. For this, taking into account Remark 3, we need to consider a lagrangian on $T^*\otimes \cal{G}$, that is an {\it action} $W=\int w(A)dx$ where $dx=dx^1\wedge ...\wedge dx^n$ and to vary it. We obtain successively:\\
 \[  \begin{array}{rcl}
 \delta A  & = & \delta a^{-1}da+a^{-1}\delta da \\
                 & = &   -(a^{-1}\delta a)(a^{-1}da)+a^{-1}d(aa^{-1}\delta a) \\
                 & = &  d(a^{-1}\delta a)+(a^{-1}da)(a^{-1}\delta a)-(a^{-1}\delta a)(a^{-1}da)
     \end{array}   \]
Setting $a^{-1}\delta a=\lambda\in {\cal{G}}={\wedge}^0T^*\otimes\cal{G}$, we thus obtain [13,14,15]:\\
\[      \delta A=d\lambda - [A,\lambda]        \]
Finally, setting $\partial w/\partial A={\cal{A}}=({\cal{A}}^i_{\tau})\in {\wedge}^{n-1}T^*\otimes \cal{G}$, we get:\\
\[  \delta W=\int {\cal{A}}\delta Adx=\int {\cal{A}}(d\lambda-[A,\lambda])dx  \]
and therefore, after integration by part, the Euler-Lagrange (EL) equations (with no right members) [1(13),15]:\\
\[     {\partial}_i{\cal{A}}^i_{\tau}+c^{\sigma}_{\rho\tau}A^{\rho}_i{\cal{A}}^i_{\sigma}=0    \]
We notice that such a linear operatorfor $\cal{A}$ has non constant coefficients linearly depending on $A$.\\
However, setting $\delta aa^{-1}=\mu$, we also have:\\
\[  \begin{array}{rcl}
\delta A  & = & -a^{-1}\delta aa^{-1}da+a^{-1}d((\delta aa^{-1})a)   \\
                & = & -a^{-1}\delta aa^{-1}da +a^{-1}d(\delta aa^{-1})a+a^{-1}\delta aa^{-1}da  \\
                & = &  Ad(a)d\mu
    \end{array}    \]
Therefore, introducing by duality $\cal{B}$ such that ${\cal{B}}\mu= {\cal{A}}\lambda$, we get the equivalent form [1 (14),15]:  \\
\[   {\partial}_i{\cal{B}}^i_{\sigma}=0   \]
which is a divergence-like operator.\\
We let the reader check by himself, as an exercise, the following other formulas for right invariant objects:\\
\[      \delta B=d\mu +[B,\mu] =Ad(a^{-1})d\lambda    \]
where $\lambda$ is used in place of $\mu$ and vice-versa, while $a^{-1}$ is used in place of $a$ and {\it signs are changed}. Similarly, setting $A\stackrel{a\rightarrow a^{-1}}{\longrightarrow}A^{-1}=ada^{-1}=-B$ and {\it caring to the sign}, we obtain therefore $dB+[B,B]=0$. This formula will be found again later on in a different framework.\\

At the end of this section, we provide the "infinitesimal" linear version of the previous results for $a$ "close" to $e$. Then there is no difference between $A$ and $B$ or between $\lambda$ and $\mu$, and we get the {\it linear gauge sequence}:\\
\[   {\wedge}^0T^*\otimes {\cal{G}}\stackrel{d}{\rightarrow}{\wedge}^1T^*\otimes {\cal{G}}\stackrel{d}{\rightarrow}{\wedge}^2T^*\otimes {\cal{G} }     \]
which is just the tensor product by $\cal{G}$ of a part of the Poincar\'e sequence for the exterior derivative.\\

\noindent
{\bf DYNAMICS ON LIE PSEUDOGROUPS}:\\

   Before using specific notations for this particular section in a way coherent with the first section, let us provide an elementary introduction to the local theory of Lie pseudogroups with notations that will also be used in the third section where they will also be coherent with the notations of the first section too.\\
   
\noindent
{\bf DEFINITION 1}: A {\it Lie group of transformations} of a manifold $X$ is a lie group $G$ with an {\it action} of $G$ on $X$ better defined by its graph $X\times G \rightarrow X\times X: (x,a)\rightarrow (x,y=ax=f(x,a))$ with the properties that $a(bx)=(ab)x$ and $ex=x, \forall x\in X, \forall a,b\in G$. \\
It is sometimes useful to distinguish the {\it source} $x$ from the {\it target} $y$ by introducing a copy $Y$ of $X$ with local coordinates $y=(y^1,...,y^n)$. Such groups of transformations have first been studied by S. Lie in 1880. Among basic examples when $n=1$ we may quote the {\it affine group} $y=ax+b$ and the {\it projective group} $y=(ax+b)/(cx+d)$ of transformations of the real line. When $n=3$ we may quote the {\it group of rigid motions} $y=ax+b$ where now $a$ is an orthogonal $3\times 3$ matrix and $b$ is a vector. Such a group is known to preserve the {\it euclidean metric} $\omega=({\omega}^{ij}={\omega}^{ji})$ and thus the quadratic form $ds^2=(dx^1)^2+(dx^2)^2+(dx^3)^2={\omega}^{ij}dx^idx^j$. When $n=4$ we may quote the {\it conformal group} of space-time with $15$ parameters (4 translations, 6 rotations, 1 dilatation, 4 elations) preserving the Minkowski metric $\omega$ or the quadratic form $ds^2=(dx^1)^2+(dx^2)^2+(dx^3)^2-c^2(dt)^2$ up to a function factor, where now $c$ is the speed of light and $t$ the time. Among the subgroups, we may consider the {\it Weyl group} with $11$ parameters preserving $\omega$ up to a constant factor and the {\it Poincar\'e group} with $10$ parameters preserving $\omega$. We recall that the three Lorentz transformations must be considered as space-time rotations, the three other rotations being pure space rotations.\\

Only ten years later, in 1890, S. Lie discovered that the Lie groups of transformations were only examples of  a wider class of groups of transformations, first called {\it infinite groups} but now called {\it Lie pseudogroups}.\\

\noindent
{\bf DEFINITION 2}: A Lie pseudogroup $\Gamma$ of transformations of a manifold $X$ is a group of transformations $y=f(x)$ solutions of a (in general nonlinear) system of partial differential (PD) equations, also called system of {\it finite Lie equations}.

Roughly, this definition means that, if ${\Pi}_q(X,Y)$ with local coordinates $(x,y_q)=(x,y^k,y^k_i,y^k_{ij},...)$ satisfying $ det(y^k_i)\neq 0$ is the q-jet bundle of invertible transformations (just replace derivatives by symbols !), there is a {\it system} ${\cal{R}}_q\subset {\Pi}_q(X,Y)$ defined by equations ${\Phi}^{\tau}(x,y_q)=0$ such that, if we have two solutions that can be composed, the composition is again a solution. However, such a point of view cannot be tested in actual practice. Instead, the idea is to use {\it sections} of ${\cal{R}}_q$, that is maps $f_q:(x)\rightarrow (x,f^k(x),f^k_i(x),f^k_{ij}(x),...)$ satisfying $det(f^k_i(x))\neq 0$ such that ${\Phi}^{\tau}(x,f_q(x)\equiv 0, \forall x\in X$. Introducing the specific section $j_q(f):(x)\rightarrow (x,f^k(x),{\partial}_if^k(x),{\partial}_{ij}f^k(x),...)$ and the composition $j_q(g\circ f)=j_q(g)\circ j_q(f)$ for maps that can be composed, an equivalent definition that can be tested is that, whenever $f_q$ and $g_q$ are sections of ${\cal{R}}_q$ that can be composed, then $g_q\circ f_q$ is also a section. A similar comment can be done for defining the inverse $f_q^{-1}$ and is left to the reader.\\

Setting now $y=x+t{\xi}(x)+...$ and passing to the limit for $t\rightarrow 0$, that is to say linearizing ${\cal{R}}_q$ around the q-jet of the identity $y=x$, we get a linear system $R_q\subset J_q(T)$ for vector fields with solutions $\Theta\subset T$ satisfying $[\Theta,\Theta]\subset \Theta$. It can be proved, for the same testing type reasons, that such a system may be endowed with a Lie algebra bracket {\it on sections} ${\xi}_q:(x)\rightarrow (x,{\xi}^k(x), {\xi}^k_i(x),{\xi}^k_{ij}(x),...)$ that we shall quickly define (see [14,15] for more details and compare to [12]). Such a bracket on sections transforms $R_q$ into a {\it Lie algebroid} and we have $[R_q,R_q]\subset R_q$. Let us first define by bilinearity $\{j_{q+1}(\xi),j_{q+1}(\eta)\}=j_q([\xi,\eta]), \forall \xi,\eta\in T$. Introducing the {\it Spencer operator} $D:R_{q+1}\rightarrow T^*\otimes R_q:{\xi}_{q+1}\rightarrow j_1({\xi}_q)-{\xi}_{q+1}$ with local components $({\partial}_i{\xi}^k-{\xi}^k_i,{\partial}_i{\xi}^k_j-{\xi}^k_{ij},...)$, we obtain the following general formula at order $q$:\\
\[   [{\xi}_q,{\eta}_q]=\{{\xi}_{q+1},{\eta}_{q+1}\}+i(\xi)D{\eta}_{q+1}-i(\eta)D{\xi}_{q+1}  ,\forall {\xi}_q,{\eta}_q\in R_q\]
where $i( )$ is the interior multiplication of a 1-form by a vector, and we let the reader check that such a definition no longer depends on the "lifts" ${\xi}_{q+1},{\eta}_{q+1}$ over ${\xi}_q,{\eta}_q$.\\

\noindent
{\bf EXAMPLE 1}: (Affine transformations)  $n=1,q=2, X={\mathbb{R}}^3$\\
With evident notations, the system ${\cal{R}}_2$ is defined by the single linear OD equation $y_{xx}=0$ and the sections are defined by $f_{xx}(x)=0$. Similarly, the sections of $R_2$ are defined by ${\xi}_{xx}(x)=0$. Accordingly, the components of $[{\xi}_2,{\eta}_2]$ at order zero, one and two are defined by the totally unusual successive formulas:\\
\[    [\xi,\eta]=\xi{\partial}_x\eta-\eta{\partial}_x\xi     \]
\[    ([{\xi}_1,{\eta}_1])_x=\xi{\partial}_x{\eta}_x-\eta{\partial}_x{\xi}_x    \]
\[    ([{\xi}_2,{\eta}_2])_{xx}={\xi}_x{\eta}_{xx}-{\eta}_x{\xi}_{xx}+\xi{\partial}_x{\eta}_{xx}-\eta{\partial}_x{\xi}_{xx}   \]
It follows that ${\xi}_{xx}=0,{\eta}_{xx}=0\Rightarrow ([{\xi}_2,{\eta}_2])_{xx}=0$ and thus $[R_2,R_2]\subset R_2$.\\

\noindent
{\bf EXAMPLE 2}: (Projective transformations)  $n=1, q=3, X={\mathbb{R}}^3$\\
The system ${\cal{R}}_3$ is defined by the single nonlinear OD equation $(y_{xxx}/y_x)-\frac{3}{2}{(y_{xx}/y_x)}^2=0$ and the sections of $R_3$ are defined by ${\xi}_{xxx}(x)=0$. The formulas for the bracket of Lie algebroid $[R_3,R_3]\subset R_3$ can be derived similarly but involve many more terms.\\

\noindent
{\bf EXAMPLE 3}: (Volume preserving transformations)  $n$ arbitrary, $q=1, X={\mathbb{R}}^n$\\
The system ${\cal{R}}_1$ is defined by the single nonlinear PD equation ${\partial}(y^1,...,y^n)/{\partial}(x^1,...,x^n)=det(y^k_i)=1$ and the sections of $R_1$ are defined by the single relation ${\xi}^i_i=0$. Accordingly, we obtain:\\
\[   ([{\xi}_1,{\eta}_1])^k_i={\xi}^r_i{\eta}^k_r-{\eta}^r_i{\xi}^k_r+{\xi}^r{\partial}_r{\eta}^k_i-{\eta}^r{\partial}_r{\xi}^k_i       \]
When summing on $k$ and $i$, the first two terms disappear (as in Example 1 !) and we get therefore $[R_1,R_1]\subset R_1$. We invite the reader to compare this result with the usual way on solutions where one defines $\Theta$ as the kernel of the Lie derivative ${\cal{L}}(\xi)\omega$ of the (volume) n-form $\omega=dx^1\wedge ...\wedge dx^n$ with respect to $\xi$ and then uses the well known formula $[{\cal{L}}(\xi),{\cal{L}}(\eta)]\omega={\cal{L}}([\xi,\eta])\omega , \forall \xi,\eta\in T$ in order to obtain $[\Theta,\Theta]\subset \Theta$.\\

\noindent
{\bf REMARK 4}: In view of the results of the first section, it is tempting to consider a Lie pseudogroup of transformations as an infinite Lie group of transformations. Such a point of view, frequently adopted in the past by russian mathematicians (V. Arnold, L.V. Ovsiannikov,...) is at the opposite of point of view adopted by the american school (D.C. Spencer and coworkers). The origin of this confusion lies in the fact that the transformations of a few Lie pseudogroups just depend on a certain number of arbitrary functions. Expanding these functions into Taylor series close to the value they have at the identity transformation allows to consider the coefficients of the series as parameters on an infinite Lie group. {\it  The best way is by far to obtain results from the finite number of defining finite or infinitesimal Lie equations, using their sections without using their solutions}. The following tricky example will illustrate this remark.\\

\noindent
{\bf EXAMPLE 4}: Let us consider the Lie pseudogroup of transformations of ${\mathbb{R}}^2$ only depending on an arbitrary function of a single variable with nonzero derivative:\\
\[      \Gamma=\{y^1=f(x^1), y^2=x^2/\frac{\partial f(x^1)}{\partial x^1}\}     \]
It is easy to check that $\Gamma$ is the set of solutions of the (involutive) nonlinear system:\\
\[   {\cal{R}}_1      \hspace{3cm}  y^1_2=0, y^2y^1_1=x^2 \Rightarrow \frac{\partial (y^1,y^2)}{\partial (x^1,x^2)} =1  \]
with corresponding system of infinitesimal Lie equations:\\
\[R_1\subset J_1(T) \hspace{3cm}  {\xi}^1_2=0, x^2{\xi}^1_1+{\xi}^2=0\Rightarrow {\xi}^1_1+{\xi}^2_2=0\]
We may thus set $f(x)=a_0+a_1x+...$ with $a_1\neq 0$ and $f$ must be the identity $f(x)=x$ whenever $y^1=x^1,y^2=x^2$.\\
Nevertheless, changing slightly the last system to the following one:\\
\[ {\cal{R}}'_1   \hspace{3cm}  y^1y^2_2-y^2y^1_2=x^1, y^1y^2_1-y^2y^1_1=-x^2\Rightarrow \frac{\partial (y^1,y^2)}{\partial (x^1,x^2)} =1  \]
does not allow now to have any generic explicit solution. We finally notice that, in the first case one can check the Lie pseudogroup property by composing explicit solutions while in the second case the corresponding Lie pseudogroup is the one preserving the {\it geometric object} $\omega=(\alpha,\beta)$ where $\alpha=x^1dx^2-x^2dx^1$ is a 1-form and $\beta=dx^1\wedge dx^2$ is a 2-form satisfying the {\it integrability condition} $d\alpha=2\beta$.\\

From now on in this section, let $X$ be the manifold of parameters as in the first section, let $Y$ be another manifold of dimension $m$ with coordinates $y=(y^1,...,y^m)$ and indices $k,l$ while $Z$ is a copy of $Y$ with coordinates $z=(z^1,...,z^m)$. Let now $\Gamma={\Gamma}_Y\subseteq aut(Y)$ be a {\it Lie pseudogroup} of transformations of $Y$ ({\it not} of $X$) of the form $z=f(y)$ with inverse $y=g(z)$. With a slight abuse of language, a map $f:X\rightarrow \Gamma : x\rightarrow f_x(y)=f(y,x)$ will be identified with a section of $X\times \Gamma$. We denote by $R_q=R_q(Y)\subseteq J_q(T(Y))$ the corresponding system of infinitesimal Lie equations. The operator ${\cal{D}}={\cal{D}}_Y:T(Y)\rightarrow J_q(T(Y))/R_q(Y)$ is a {\it Lie operator}, that is an operator such that ${\cal{D}}{\xi}_1=0, {\cal{D}}{\xi}_2=0\Rightarrow {\cal{D}}[{\xi}_1,{\xi}_2]=0$, where the bracket is the standard bracket of vector fields on $Y$. Again with a slight abuse of language, if $\Theta=\Theta_Y\subseteq T(Y)$ is the sheaf of germs of solutions of $\cal{D}$, we have $[\Theta,\Theta]\subset \Theta$.\\

The main problem that we have now to solve is to construct an analogue of the gauge sequence with {\it left} and {\it right} points of view, both with a corresponding variational calculus. For this, considering $f_{x+dx}\circ f_x^{-1}-id$ and passing to the limit, we obtain a vector valued 1-form (CARE) $v\in T^*\otimes {\Theta}_Z\subset T^*\otimes V(X\times Z)$ called {\it generalized speed} and defined by the formula:\\
\[       v(z,x)=\frac{\partial f}{\partial x}(g(z,x),x)    \]
with no one of the indices for simplicity. Such a definition is just similar to the definition of the eulerian speed $v(x,t)$ in continuum mechanics for a movement $x=f(x_0,t)$ of a point/particle with position $x_0$ at time $t_0$, $x$ at time $t$ and $x+v(x,t)dt$ at time $t+dt$. However, the Lie pseudogroup background is not so evident and the reader may now, as a motivation, get in mind the composition $z\rightarrow y\rightarrow z+v(z,x)dx$.\\
Similarly, considering $f_x^{-1}\circ f_{x+dx}-id$ and passing to the limit, we may obtain the pull-back $u\in T^*\otimes {\Theta}_Y\subset T^*\otimes V(X\times Y)$ defined by the formula:\\
\[       \frac{\partial f}{\partial x}(y,x)\equiv v(f(y,x),x)\equiv \frac {\partial f}{\partial y}(y,x)u(y,x)   \]
again with no one of the indices for simplicity as only the tangent mapping is involved.\\
Conversely, if $y=g(z,x),\bar{y}=\bar{g}(z,x)\in \Gamma$, then, from the implicit function theorem, we get $\bar{y}=h(y,x)\in \Gamma$ in general. However, from the identities:\\
\[ f(g((z,x),x)\equiv z \equiv \bar{f}(\bar{g}(z,x),x) \hspace{5mm} , \hspace{5mm} \bar{f}(h(y,x),x)\equiv f(y,x)    \hspace{5mm}, \hspace{5mm}  \bar{g}(z,x)\equiv h(g(z,x),x)  \]
we deduce the relations:\\
\[ -v= \frac{\partial f}{\partial y}\frac{\partial g}{\partial x} \hspace{2cm}, \hspace{2cm}-\bar{v}=\frac{\partial\bar{f}}{\partial\bar{y}}\frac{\partial\bar{g}}{\partial x}   \]
and thus:\\
\[   v(z,x)=\bar{v}(z,x)  \Leftrightarrow \frac{\partial\bar{g}}{\partial x}=\frac{\partial h}{\partial y}\frac{\partial g}{\partial x}  \Leftrightarrow \frac{\partial h}{\partial x}=0  \Leftrightarrow \bar{y}=h(y)\in \Gamma  \]
a result showing the right invariance of the generalized speed.\\
Exactly as in the previous section, we notice that the passage from $v$ to $u$ just amounts to change $f$ to $g$ while changing the sign (CARE) as we have indeed the identity:\\
\[    y\equiv g(f(y,x),x)  , \forall x\in X   \]
and, differentiating with respect to $x$, we get:\\
\[  \frac{\partial g}{\partial z}(z,x)v(z,x)+\frac{\partial g}{\partial x}(z,x)\equiv 0    \]
and thus:\\
\[     u(g(z,x),x)\equiv -\frac{\partial g}{\partial x}(z,x)  , \forall x\in X   \]
or equivalently:\\
\[     u(y,x)\equiv -\frac{\partial g}{\partial x}(f(y,x),x)  , \forall x\in X   \]
but $u,v\in T^*\otimes \Theta$ do not provide the same sections.\\

The next important step will be to provide the {\it compatibility conditions}. For this, introducing:\\
\[     v^k_i=\frac{\partial f^k}{\partial x^i}(g(z,x),x)   \]
we successively obtain from the chain rule for derivatives and the above comment:\\
\[   \begin{array}{rcl}
   \frac{\partial v^k_i}{\partial x^j} & = & \frac{{\partial}^2f^k}{\partial x^i\partial y^l}\frac{\partial g^l}{\partial x^j}+\frac{{\partial}^2f^k}{\partial x^i\partial x^j}  \\
   &    &   \\
     & = & \frac{\partial v^k_i}{\partial z^u}\frac{\partial f^u}{\partial y^l}\frac{\partial g^l}{\partial x^j}+\frac{{\partial}^2 f^k}{\partial x^i\partial x^j}  \\
       &  &   \\
      &  =  & -v^u_j\frac{\partial v^k_i}{\partial z^u} + \frac{{\partial}^2f^k}{\partial x^i\partial x^j}
      \end{array}  \]
Exchanging $i$ and $j$, then substracting, we get the formula:\\
\[     {\partial}_iv_j -{\partial}_jv_i+[v_i,v_j]=0    \]
that we may rewrite in the condensed form $dv+[v,v]=0$ as $v$ can be considered as an horizontal form on  $X\times Z$. Also, as in the first section, with $u$ in place of $A$, $v$ in place of $B$ and $\frac{\partial g}{\partial x}$ in place of $B^{-1}=-A$, while caring about the sign, we have similarly $du-[u,u]=0$.\\
Collecting these results and taking into account that $\cal{D}$ is a Lie operator with ${\partial}_i{\cal{D}}=0, \forall i=1,...,n$, we obtain therefore:\\

\noindent
{\bf THEOREM 2}: There exists a differential sequence:\\
\[ X\times \Gamma \longrightarrow T^*\otimes \Theta \longrightarrow {\wedge}^2T^*\otimes \Theta  \]

The hard step in this paper will be to construct the corresponding variational calculus with constraint. For this, $\delta$ being the usual symbol for variation, we may introduce the vertical vector field $\eta=\eta(z,x)\in \Theta$ such that $\eta (f(y,x),x)=\delta f(y,x)$ and we get (compare to [7,14,20]):\\
\[  \delta v+\eta \frac{\partial v}{\partial z}=\frac{\partial \eta}{\partial x}+v\frac{\partial \eta}{\partial z} =\frac{d\eta}{dx}
\hspace{1cm} \Rightarrow \hspace{1cm}    \delta v=\frac{\partial \eta}{\partial x}+[v,\eta]   \]
Exactly as above and for a later use, by analogy with continuum mechanics, we may introduce the Jacobian determinant:\\
\[  \Delta (z,x)=\frac{\partial (f)}{\partial (y)}(g(z,x),x)=det(\frac{\partial f}{\partial y})(g(z,x),x) \]
Using the well known formulas for jacobian determinants:\\
\[ \frac{\partial (f+df)}{\partial (y)}/\frac{\partial (f)}{\partial (y)}=\frac{\partial (z+vdx)}{\partial (z)}=1+\frac{\partial v}{\partial z}dx+...   \] 
and passing to the limit, we obtain therefore $\frac{d\Delta}{dx}=\Delta \frac{\partial v}{\partial z}$ and thus:\\
\[ \delta \Delta+\eta \frac{\partial \Delta}{\partial z}=\Delta \frac{\partial \eta}{\partial z} \Rightarrow \delta \Delta =\Delta \frac{\partial \eta}{\partial z}-\eta \frac{\partial \Delta}{\partial z}   \]
In order to overcome the minus sign in this formula and simplify it, we may therefore introduce the analogue $\rho=1/\Delta$ of the mass per unit of volume and obtain the variation [5,12,18]:\\
\[  \delta \rho +\eta\frac{\partial \rho}{\partial z}=-\rho \frac{\partial \eta}{\partial z} \hspace{1cm}\Rightarrow \hspace{1cm} \delta \rho=-\frac{\partial (\rho {\eta}^k)}{\partial z^k} . \]
We notice that such a result could not be obtained in [1] where $\Delta=1$ {\it by assumption}, though this is coherent with the divergence-free condition for the speed in this case.\\

As in the first section, we may also introduce the pull back of $\eta$ as the vertical vector $\xi=\xi(y,x)\in \Theta$ defined by the formula:\\
\[  \eta (f(y,x),x)\equiv \frac{\partial f}{\partial y}(y,x)\xi (y,x)  \hspace{5mm}\Leftrightarrow\hspace{5mm}
\frac{\partial g}{\partial z}(z,x)\eta (z,x)\equiv \xi (g(z,x),x)  \]
exactly as we got $u$ from $v$ but now with $\delta$ in place of $d$. We obain the following {\it tricky theorem}, coherent with the results of the first section but absent from [1,2]:\\

\noindent
{\bf THEOREM 3}:     $\hspace{2cm}  \delta v=\frac{\partial f}{\partial y} \frac{\partial \xi}{\partial x} $.\\

\noindent
{\it Proof}: Using the well known result saying that the bracket of vector fields commutes with the action of any diffeomorphism while working ONLY with $(z,x)$ and CARING ABOUT THE SIGN, we successively get from the above formulas:\\
\[  \begin{array}{rcl}
\frac{\partial g}{\partial z}(z,x)\delta v & = & \frac{\partial g}{\partial z}\frac{\partial \eta}{\partial x}+\frac{\partial g}{\partial z} [v,\eta]  \\
  &  &  \\
  & = & \frac{\partial g}{\partial z}\frac{\partial \eta}{\partial x}+[\frac{\partial g}{\partial z} v,\frac{\partial g}{\partial z} \eta ]  \\
    &  &  \\
    &  =  & \frac{\partial g}{\partial z} \frac{\partial \eta}{\partial x}-[\frac{\partial g}{\partial x},\xi]  \\
      &  &  \\
    &  = &  \frac{\partial}{\partial x}(\frac{\partial g}{\partial z}\eta)-\frac{\partial}{\partial x}(\frac{\partial g}{\partial z})\eta-[\frac{\partial g}{\partial x},\xi]  \\
    &  &  \\
    &  = & \frac{\partial \xi}{\partial x}+\frac{\partial \xi}{\partial y}\frac{\partial g}{\partial x} - \frac{\partial}{\partial z}(\frac{\partial g}{\partial x})\eta -\frac{\partial g}{\partial x}\frac{\partial \xi}{\partial y}+\eta\frac{\partial}{\partial z}(\frac{\partial g}{\partial x})  \\
      &  &  \\
    &  = & \frac{\partial \xi}{\partial x}(g(z,x),x)   
    \end{array}   \]
as we have indeed:\\
\[  \xi \frac{\partial}{\partial y}=\eta \frac{\partial}{\partial z} \hspace{1cm},\hspace{1cm} \frac{d}{dx}\xi (g(z,x),x)\equiv \frac{\partial \xi}{\partial x}+\frac{\partial \xi}{\partial y}\frac{\partial g}{\partial x}   \]
thus ending the proof.\\
\hspace*{12cm}     Q.E.D.  \\

By analogy with continuum mechanics, among general action integrals of the form $W=\int w(\frac{\partial f}{\partial x},\frac{\partial f}{\partial y})dy\wedge dx$, we shall only consider the action integrals of the form $W=\int \rho (z,x)w(v(z,x))dz\wedge dx$ and vary them. As in [7(V),14], the best procedure in actual practice is to introduce a "small" parameter $\epsilon$ and consider the family $f(y,x;\epsilon)\leftrightarrow g(z,x;\epsilon)$ with variation $\delta=\partial /\partial\epsilon$ also leading to the previous results when $(z,x;\epsilon)$ are considered as independent variables. Of course, we have to add the constraint $z=f_x(y)=f(y,x)\in \Gamma$ and this will be the delicate point to overcome.\\
Introducing ${\cal{V}}=\frac{\partial w}{\partial v}=({\cal{V}}^i_k)$ or its dual $\cal{U}$ we get:\\
\[     \delta (\rho w)=-\frac{\partial}{\partial {z}^k}(\rho w{\eta}^k)+\rho {\cal{V}}^i_k\frac{\partial {\eta}^k}{\partial {x}^i}+\rho {\cal{V}}^i_kv^l_i\frac{\partial {\eta}^k}{\partial {z}^l}.  \]
Using $x$ in place of $\epsilon$ for $\rho$, we get the "{\it conservation of mass}" identity:\\
\[    \frac{\partial \rho}{\partial {x}^i}+\frac{\partial}{\partial {z}^k}(\rho {v}^k_i)\equiv 0  \]
Finally, integrating by part over $X\times Z$ exactly as for Eulerian coordinates $(x,t)$ in classical hydrodynamics [7,14,20], we obtain the EL equations over $X\times Z$, up to sign:\\
\[    {\cal{E}}\eta \equiv   \rho (\frac{\partial {\cal{V}}^i_k}{\partial {x}^i}+{v}^l_i\frac{\partial {\cal{V}}^i_k}{\partial {z}^l}){\eta}^k=0, \hspace{1cm} \forall \eta \in \Theta\]
or equivalently over $X\times Y$:\\
\[    \frac{\partial {\cal{U}}^i_k}{\partial {x}^i}{\xi}^k=0, \hspace{1cm} \forall \xi \in \Theta.  \]
It is just the solution of this problem that we have already solved for creating PD {\it optimal control theory}. We consider two situations:\\
 
1)  If $\cal{D}$ can be parametrized by an operator ${\cal{D}}_{-1}$, then any $\eta\in\Theta$ can be written in the form ${\cal{D}}_{-1}\theta=\eta$ for an arbitrary $\theta\in V(X\times Z)$ in the previous variation. Then, denoting by $ad$ the {\it formal adjoint} of an operator and integrating by part with respect to $z$, we find $ad({\cal{D}}_{-1}){\cal{E}}=0$ as final EL equations.\\

2)  In the general situation, if $\{{\Phi}^{\tau}\}$ is a fundamental generating set of differential invariants of $\Gamma$ at order $q$ and $\{{\omega}^{\tau}(y)\}$ their value at the identity, the {\it Lie form} of the {\it finite Lie equations} defining $\Gamma$ at order $q$ is ${\Phi}^{\tau}-{\omega}^{\tau}(y)=0$ in the sense that we have a defining {\it groupoid} of order $q$. It follows that the {\it constrained} variational calculus for an action density $w$ amounts to vary the {\it unconstrained} action density $w+{\lambda}_{\tau}({\Phi}^{\tau}-{\omega}^{\tau}(y))$ where $\{{\lambda}_{\tau}\}$ are {\it Lagrange multipliers}. As the vertical bundle of the groupoid is isomorphic to the corresponding infinitesimal Lie equations over the target, changing $\lambda$ to $\bar{\lambda}$ by duality, the variation becomes $\delta w+\bar{\lambda}{\cal{D}}\eta$. Integrating by part as usual while caring about the signs, we obtain the EL equations ${\cal{E}} -ad({\cal{D}})\bar{\lambda}=0$. Finally, we just need to eliminate the Lagrange multipliers by introducing compatibilty conditions for $ad({\cal{D}})$ of the form $ad({\cal{D}}_{-1})$ and obtain the same EL equations as before but with a slightly different approach.\\

\noindent
{\bf EXAMPLE 5}: When $n=1, m=3, Y={\mathbb{R}}^3$ with Euclidean metric, time $t$, initial position $x_0$, final position $x$ and transformations of the form $x=f(x_0,t)$, the above $v(x,t)$ is just the ordinary speed $\vec{v}$ in Euler variables and $u(x_0,t)$ its differential pull back in Lagrange variables. If $\Gamma=aut(Y)$ and the action density is the standard kinetic energy $w=\frac{1}{2}\rho {\vec{v}}^2$ where $\rho$ is the mass per unit volume, the EL equations become, up to sign [5(V),12]:\\
\[    \rho\vec{\gamma}\equiv \rho \frac{d\vec{v}}{dt}\equiv \rho (\frac{\partial \vec{v}}{\partial t}+(\vec{v}.\vec{\nabla})\vec{v})\equiv \rho(\frac{\partial \vec{v}}{\partial t}+(\vec{\nabla}\wedge \vec{v})\wedge \vec{v}+\vec{\nabla}(\frac{1}{2}{\vec{v}}^2))=0.  \]
Now, if $\Gamma$ is the Lie pseudogroup of volume preserving transformations defined by $\Delta=1\Leftrightarrow \rho=1$, then $\Theta$ is the Lie algebra of divergence-free vector fields and we get:\\
\[    \vec{\gamma}.\vec{\eta}=0, \hspace{1cm} \forall \eta\in\Theta \Leftrightarrow \vec{\nabla}.\vec{\eta}=0\]
In {\it that specific case}, ${\cal{D}}=div$ can be parametrized by the $curl$ operator, that is $\vec{\eta}=\vec{\nabla}\wedge \vec{\theta} $ for an arbitrary vector field $\vec{\theta}$. As the $curl $ operator is self-adjoint up to sign, introducing the so-called "{\it vortex}" vector $\vec{\omega}=\frac{1}{2}\vec{\nabla}\wedge \vec{v}$ (see below for the numerical factor), we get at once the EL equations of Arnold, namely:\\
\[     \frac{1}{2}\vec{\nabla}\wedge \vec{\gamma}\equiv \frac{\partial \vec{\omega}}{\partial t}+[\vec{v},\vec{\omega}]=0  \]
in a {\it purely formal way}.\\
We point out that the numerical factor 1/2 comes from the fact that, in the case of a rigid body motion $x=a(t)x_0+b(t)$, then $v=\dot{x}=\dot{a}x_0+\dot{b}=\dot{a}a^{-1}x+(\dot{b}-\dot{a}a^{-1}b)$. An easy computation in local coordinates then shows that $\dot{a}a^{-1}$ is a skew-symmetric $3\times 3$ matrix amounting to the vector product by the {\it vortex vector} $\vec{\omega}$ as previously defined. \\Otherwise, considering the action density $w+\lambda (\Delta -1)$ and varying it, we get $\delta w+\lambda \vec{\nabla}.\vec{\eta}$. Integrating by part, we get $\vec{\gamma}+\vec{\nabla}\lambda=0\Rightarrow \vec{\nabla}\wedge\vec{\gamma}=0$ as before but the Lagrange multiplier $\bar{\lambda}=\Delta\lambda=\lambda$ now plays the part of the {\it pressure}, as well known in fluid dynamics.\\

\noindent
{\bf EXAMPLE 6}: With $m=1$, setting $z'=\partial z/\partial y$ and so on, Le Lie form of the group of projective transformations of the real line is $\Phi\equiv \frac{z'''}{z'}-\frac{3}{2}(\frac{z"}{z'})^2=0$ where $\Phi$ is the well known Schwarzian derivative. Using the prolongation formulas $\delta z=\eta, \delta z'={\eta}'z'=(\delta z)',...$ and so on, it is not so easy to check that $ \delta \Phi=(z')^2{\eta}'''$ and the needed isomorphism is thus $\bar{\lambda}=(z')^2\lambda$.\\

\noindent
{\bf EXAMPLE 7}: From everybody diving experience or looking at corks on the surface of the sea, a natural model for swell is well described by water particles moving along circles with radius smaller and smaller with depth. A parametrized form in the vertical $(x,y)$ plane can be:\\
\[    x=a+R(b)cos(\omega t-\varphi (a)) \hspace{1cm},\hspace{1cm}y=b+R(b)sin(\omega t-\varphi (a))  \]
As we have the well known formula $\frac{\partial (x,y)}{\partial (x_0,y_0)}=\frac{\partial (x,y)}{\partial (a,b)}/\frac{\partial (x_0,y_0)}{\partial (a,b)}$, the movement is that of an uncompressible fluid iff $\frac{\partial (x,y)}{\partial (a,b)}=\frac{\partial (x_0,y_0)}{\partial (a,b)}=1+RR'{\varphi}'+(R{\varphi}'+R')sin$ does not depend on time, that is to say ${\varphi}'(a)=-R'/R(b)=k$ and thus $\varphi(a)=ka+c, R(b)=R_0e^{-kb}$. In that case, it is {\it not at all evident to see directly} that the movement is NOT stationary in the fixed frame. However, for proving this result, one just needs to notice that the {\it trajectories} are bounded circles while the {\it stream lines} are infinite cycloidal curves (consider the surface of the sea !). Finally, the movement becomes stationary in the moving frame $\bar{x}=x-\frac{\omega}{k}t,\bar{y}=y$ associated with a boat "surfing" on the crest of a given wave. Indeed, changing the parameters $(a,b)$ 
to $(\bar{a}=a-\frac{\omega}{k}t,\bar{b}=b)$, we notice that the new parametrization $(\bar{a},\bar{b})\rightarrow (\bar{x},\bar{y})$ does not involve time $t$ anymore.\\

\noindent
{\bf DYNAMICS ON LIE GROUPOIDS}:\\

Let us revisit the variation formula provided by Theorem 3, changing the notations when $n=1$, one parameter $t$ and movement $x=f(x_0,t)$ with inverse $x_0=g(x,t)$ as usual in continum mechanics where the speed in Eulerian coordinates $(x,t)$ is defined by :\\
\[    v(f(x_0,t),t)\equiv \frac{\partial f}{\partial t}(x_0,t) .  \]
Using space-time coordinates $(x^1,x^2,x^3,x^4)$ with $x^4=ct$ and reference speed of light $c=1$, we may divide $(dx^1,dx^2,dx^3,dx^4)$ by $dt$ and obtain the {\it extended speed} $(v^1,v^2,v^3, 1)=
(\vec{v},1)$ in Special Relativity with variation $(\delta \vec{v},0)$. Introducing the "time" variable $t_0$ corresponding to $x_0$ and the function $t_0=g^4(x,t)$, we may finally consider the transformation $x=f(x_0)$ with inverse $x_0=g(x)$, but now on space-time, and try to obtain $\delta v=[v,\eta]$ on space-time, {\it on the condition to get rid of} $t_0$. For this we refer to [8,12] as it is out of the scope of this paper. In particular, the strange form $dx\wedge dt$ is no longer a wedge product of forms but a volume form on space-time as we have no longer TWO separate manifolds but only ONE, a result leading to relativistic mechanics. \\
Having in mind the first section where $G$ was an abstract group, we have now a Lie pseudogroup of transformations of a manifold, a particularly simple example being produced by a Lie group of transformations defined by the graph  $X\times G\rightarrow X\times X$ of an action. In this case, it is known [14,15] that the (nonlinear) Spencer sequence is isomorphic to the gauge sequence and the Spencer (nonlinear) sequence is therefore the ONLY candidate for a generalization to cases where no Lie group action is involved, like the Lie pseudogroup of volume preserving, symplectic or contact transformations.\\
We have already explained in many books [14,15] that  such an aproach, where sections of Lie groupoids/jet bundles generalize the gauging approach, is JUST the way to understand the Cosserat theory of continuum mechanics and Weyl theory of electromagnetism. The only choice remaining on space-time is to select a convenient "{\it group candidate}"  and we refer to the previous books for discovering why it MUST be the conformal group of space-time endowed with the Minkowski metric.\\
Finally, exactly like in the first section, a linearized version exists and we refer the reader to [16,17,19] for a simple presentation of this framework with applications to field-matter coupling.\\

 Coming back to mathematics for a few lines, there is another differential sequence to be found in the literature and that we did not speak about, namely the {\it Janet sequence}, though it is for sure the best known differential sequence. For short, if $E,F, F_0,F_1,...$ denote vector bundles over $X$, we use the same letters for the corresponding sets (sheaves to be exact) of sections and such an interpretation must be used whenever operators are involved. Starting from a vector bundle $E$ (for example $T$) and a linear differential operator ${\cal{D}}:E\rightarrow F:\xi\rightarrow \eta$, if we want to solve the linear system with second member ${\cal{D}}\xi=\eta$ even locally, one needs "{\it compatibility conditions}" in the form ${\cal{D}}_1\eta=0$. Denoting now $F$ by $F_0$, we may therefore look for an operator ${\cal{D}}_1:F_0\rightarrow F_1:\eta\rightarrow \zeta$ and so on. Under assumptions on $\cal{D}$ which are out of the scope of this paper (involutivity !), the french mathematician M. Janet has proved in 1920 that such a chain of operators ends after $n$ steps and we obtain the {\it linear Janet sequence}, namely [15]:\\
 \[   0\rightarrow \Theta\rightarrow E\stackrel{\cal{D}}{\rightarrow}F_0\stackrel{{\cal{D}}_1}{\rightarrow}F_1\stackrel{{\cal{D}}_2}{\rightarrow}...\stackrel{{\cal{D}}_n}{\rightarrow} F_n\rightarrow 0  \]
 It follows that we only have at our disposal for any application where group theory seems to be involved, three linear differential sequences, namely the Janet sequence, the Spencer sequence and the gauge sequence. As these sequences are made by quite different operators, {\it the use of one excludes the use of the others}.\\
 
 In order to escape from this dilemna, at the end of this paper and for the sake of clarifying the key idea of the brothers Cosserat by using these new mathematical tools, we shall explain, in a way as elementary as possible while using only the linear framework, why THE JANET SEQUENCE AND THE GAUGE SEQUENCE CANNOT BE USED IN CONTINUUM MECHANICS. By this way we hope to convince the reader about the need to use another differential sequence, namely the SPENCER SEQUENCE, though striking it could be. Also we shall use very illuminating examples in order to illustrate our comments.\\
First of all we exhibit the isomoprphism existing between the linear gauge sequence and the linear Spencer sequence. For this, if {\it now} $G$ {\it acts on} $X$ with a basis ${{\xi}_{\tau}=\{\xi}^k_{\tau}{\partial}_k\} $ of infinitesimal generators, we may introduce the bundle $J_q(T)$ of $q$-jets of $T$ over $X$, that is the vector bundle over $X$ with sections transforming like the derivatives of vector fields up to order $q$, and the map:\\
\[     {\wedge}^0T^*\otimes {\cal{G}}\rightarrow J_q(T): {\lambda}^{\tau}(x) \rightarrow {\lambda}^{\tau}(x){\partial}_{\mu}{\xi}^k_{\tau}(x)     \]
It is known [14,p. 308] that this map becomes injective for $q$ large enough and we may call $R_q$ its image for such a $q$. It follows from its definition that $R_q\simeq R_{q+1}$ is a system of infinitesimal Lie equations of finite type and we get for the Spencer operator [12,15,17]:\\
\[     D:R_{q+1}\rightarrow T^*\otimes R_q:{\xi}_{q+1}\rightarrow ({\partial}_i{\xi}^k_{\mu}-{\xi}^k_{\mu+1_i})={\partial}_i{\lambda}^{\tau}(x){\partial}_{\mu}{\xi}^k_{\tau}(x)  \]
Accordingly, the linear gauge sequence is isomorphic to the linear Spencer sequence:\\
\[  0\rightarrow \Theta \rightarrow {\wedge}^0T^*\otimes R_q\stackrel{D}{\rightarrow} {\wedge}^1T^*\otimes R_q\stackrel{D}{\rightarrow} {\wedge}^2T^*\otimes R_q  \]
the three isomorphisms being induced by the (local) isomorphism $X\times {\cal{G}}\rightarrow R_q$ just described above. It is essential to notice that, though the linear Spencer sequence and the isomorphisms crucially depend on the action, by a kind of "{\it miracle}" the linear gauge sequence no longer depends on the action.\\
This result proves that {\it the linear Spencer sequence generalizes the linear gauge sequence},  with the major gain that it can be used even for Lie pseudogroups of transformations that are not coming from Lie groups of transformations, as we shall see eamples in the sequel.\\

\noindent
{\bf REMARK 5}: When $n=3$ and we deal with the Lie group of rigid motions, the corresponding EL equations are nothing else but the formal adjoint of the (first) Spencer operator. Surprisingly, this is EXACTLY the result found by the brothers Cosserat, namely the so-called {\it stress} and {\it couple-stress} equations for Cosserat media [9, p 137,14,16,18,19].\\

\noindent
{\bf REMARK 6}: The above result, in perfect agreement with the piezzoelectric or photoelastic coupling of elasticity and electromagnetism, CONTRADICTS gauge theory where the lagrangians are functions on ${\wedge}^2T^*\otimes \cal{G}$ and NOT on $T^*\otimes\cal{G}$ as in the previous remark [14,16].\\

  Let us consider a (finite) volume ${\int}_VdV$ in ${\mathbb{R}}^3$ limited by a (closed) surface $S={\int}_SdS$ and let us introduce the outside unit normal (pseudo) vector $\vec{n}=(n_j)$ on $S$. Let us now suppose that the surface element $dS$ is acted on by the outside with a force $d\vec{F}=\vec{\sigma}dS$ and a couple $d\vec{C}=\vec{\mu}dS$, where both $\vec{\sigma}$ and $\vec{\mu}$ linearly depend on $\vec{n}$ through the {\it stress} tensor density $\sigma=({\sigma}^{ij})$ and the {\it couple-stress} tensor density $\mu=({\mu}^{r,ij}=-{\mu}^{r,ji})$. It must be noticed that, using the standard Cauchy tetrahedral device, there is no reason "a -priori" to suppose that the stress tensor is symmetric. We also suppose that   the volume element $dV$ is acted on by (see later on for the sign) a force $-\vec{f}dV$ and a momentum $-\vec{m}dV$ with $\vec{f}=(f^j)$ and $\vec{m}=(m^{ij}=-m^{ji})$.\\
  Our purpose is now to study the equilibrium of the corresponding torsor fields with respect to an arbitrary cartesian frame $0x^1x^2x^3$.\\
   The equilibrium of forces is satisfied if we have the relation:\\
   \[      {\int}_S\vec{\sigma}dS-{\int}_V\vec{f}dV=0  \Rightarrow {\int}_S{\sigma}^{ij}n_idS-{\int}_Vf^jdV=0\]
Using Stokes formula, this is equivalent to the well known {\it stress equations}:\\
\[        {\partial}_i{\sigma}^{ij}=f^j   \]
This result shows that the surface density of forces $\vec{\sigma}$ is equivalent, from the point of view of force equilibrium, to a volume density of forces $\vec{f}$ and this interpretation explains the sign adopted.\\
   Finally, the equilibrium of forces being satisfied, it is known that the equilibrium of momenta is also satisfied if it is satisfied with respect to an arbitrarily chosen cartesian frame. Hence, introducing the vector $\vec{r}=(x^1,x^2,x^3)$, the equilibrium of momenta is satisfied if we have the relation:\\
   \[     {\int}_S(\vec{\mu}+\vec{r}\wedge \vec{\sigma})dS-{\int}_V(\vec{m}+\vec{r}\wedge \vec{f})dV=0 \]
   Projecting onto the axis $Ox^3$, we obtain:\\
   \[  {\int}_S({\mu}^{r,12}+x^1{\sigma}^{r2}-x^2{\sigma}^{r1})n_rdS-{\int}_V(m^{12}+x^1f^2-x^2f^1)dV=0\]
   Using again Stokes formula and the previous stress equations, we obtain the {\it couple-stress equations}:\\ 
   \[      {\partial}_r{\mu}^{r,ij}+{\sigma}^{ij}-{\sigma}^{ji}=m^{ij}  \]
 This result shows that the surface density of forces $\vec{\sigma}$ and couples $\vec{\mu}$ is equivalent, from the point of view of torsor equilibrium, to a volume density of forces $\vec{f}$ and to a volume density of momenta $\vec{m}$, provided the preceding stress and couple-stress equations are satisfied, and this interpretation explains the sign adopted.\\
The combination of the stress AND couple-stress equations have first been exhibited by E. and F. Cosserat in 1909 [8,9,p137] WITHOUT ANY STATIC EQUILIBRIUM EXPERIMENTAL BACKGROUND and we now explain the key argument leading to the same equations just from group theoretical arguments. Of course, most of the engineering continua such as steel, concrete, glass, wate,... have the specific "{\it constitutive laws}"  $\mu=0, m=0$ and we obtain therefore ${\sigma}^{ij}={\sigma}^{ji}$, that is the stress tensor is symmetric, a situation not always encountered in liquid crystals.\\

First of all, for the reader not familiar with the Spencer operator, we exhibit a similar result in a quite simpler 1-dimensional situation that will allow to recapitulate all the previous results..\\

\noindent
{\bf EXAMPLE 8}: Let us consider the Lie group of affine transformations of the real line defined by the group action $y=a^1x+a^2$. The corresponding 2-dimensional Lie group $G$ has coordinates $a=(a^1,a^2)$ and the group composition law is $ab=(a^1,a^2)(b^1,b^2)=(a^1b^1,a^1b^2+a^2)$ with inverse law $a^{-1}=(1/a^1,-a^2/a^1)$. Accordingly, we have $a^{-1}da=((1/a^1)da^1,(1/a^1)da^2)$ and obtain at once the Maurer-Cartan forms ${\omega}^1=(1/a^1)da^1, {\omega}^2=(1/a^1)da^2$ with the two Maurer-Cartan equations $d{\omega}^1=0, d{\omega}^2+{\omega}^1\wedge{\omega}^2=0$. Finally, a corresponding basis of infinitesimal generators of the action may be obtained with ${\xi}^1=x\frac{\partial}{\partial x}, {\xi}^2=\frac{\partial}{\partial x}$ and we have in a coherent way $[{\xi}^1,{\xi}^2]=-{\xi}^2$, that is the only non zero structure constant is $c^2_{12}=-1$. It follows from the first section that the resulting EL-equations are either:\\
\[  {\partial}_i{\cal{A}}^i_1+A^2_i{\cal{A}}^i_2=0\hspace{1cm},\hspace{1cm} {\partial}_i{\cal{A}}^i_2-A^1_i{\cal{A}}^i_2=0   \]
or simply:\\
\[  {\partial}_i{\cal{B}}^i_1=0 \hspace{1cm},\hspace{1cm}  {\partial}_i{\cal{B}}^i_2=0  \]
if we use a pull-back by the adjoint action. In both cases these equations could not have anything to do with the stress and couple-stress equations previously exhibited.\\
Let us now deal with the Spencer sequence instead of the gauge sequence in this framework.\\
First of all, we may consider the above Lie group of transformations as a Lie pseudogroup defined by the second order system of finite Lie equations $y_{xx}=0$. The corresponding system $R_2\subset J_2(T)$  of infinitesimal Lie equations is ${\xi}_{xx}=0$ and the isomorphisms between the gauge sequence and the Spencer sequence is induced by the maps:\\
\[  ({\lambda}^1(x),{\lambda}^2(x))\rightarrow ({\xi}(x)=x{\lambda}^1(x)+{\lambda}^2(x),{\xi}_x(x)={\lambda}^1(x),{\xi}_{xx}(x)=0)   \]
The only two non-zero components of the Spencer operator become:\\
\[  {\partial}_x{\xi}(x)-{\xi}_x(x)=x{\partial}_x{\lambda}^1(x)+{\partial}_x{\lambda}^2(x), {\partial}_x{\xi}_x(x)-0={\partial}_x{\xi}_x(x)={\partial}_x{\lambda}^1(x) \]
Equating to zero these two components amounts to have:\\
\[    {\partial}_x{\lambda}^1=0 \hspace{1cm},\hspace{1cm}  {\partial}_x{\lambda}^2=0  \]
Accordingly, gauging $\lambda$ just amounts to choose an arbitrary section of $R_2$.\\
The final touch, that could not be in the mind of any reader even on this very simple example, is to work out the formal adjoint of the Spencer operator. For this, multiplying the first component by a test function ${\sigma}(x)$, the second by a test function ${\mu}(x)$, then summing and integrating by parts, we get the EL-equations (up to sign) as the kernel of the following operator with second members $(f,m)$:\\
\[      {\partial}_x{\sigma}=f  \hspace{1cm},\hspace{1cm}  {\partial}_x\mu +\sigma=m   \]
The comparison with the previous mechanical results needs no comment.\\

Taking into account this example, we now study the foundation of elasticity theory and we restrict the study to 2-dimensional (infinitesimal) elasticity for simplicity as the general situation has already been treated elswhere and we just want to explain why {\it the only founding problem of elasticity is the choice of an underlying Lie pseudogroup and an adapted differential sequence}.\\

\noindent
1){\it The gauge sequence cannot be used}: \\
Looking at the book [9] written by E. and F. Cosserat, it seems {\it at first sight} that they just construct the first operator of the nonlinear gauge sequence for one parameter [9, p 7], two parameters [9, p 66], three parameters [9, p 123] and finally four space-time parameters [9, p 189]. This is NOT TRUE indeed because, according to the comment done in the first section or Example, the EL-equations are either a divergence like operator or a linear operator with coefficients depending on $A$, a situation not met in the couple-stress equations which is a linear operator with constant coefficients, not of divergence type. In fact, a carefull study of the book proves that somewhere the action of the group on the space is used, but this is well hidden among many very technical formulas (Compare [9] p 136 with [14] p 295).\\

\noindent
2){\it The Janet sequence cannot be used}:\\
This result is even more striking because ALL texbooks of elasticity use it along the same scheme that we now describe. Indeed, after gauging the translation by defining the "{\it displacement field}"  $\xi=({\xi}^1(x),{\xi}^2(x))$ of the body, from the initial point $x=(x^1,x^2)$ to the point $y=x+\xi(x)$, one introduces the (small) "{\it deformation tensor}" $\epsilon=1/2{\cal{L}}(\xi)\omega$ as one half the Lie derivative with respect to $\xi$ of the euclidean metric $\omega$, namely, in our case, the three components (care):\\
\[ \epsilon=({\epsilon}_{11}={\partial}_1{\xi}^1,{\epsilon}_{12}={\epsilon}_{21}=1/2({\partial}_1{\xi}^2+{\partial}_2{\xi}^1),{\epsilon}_{22}={\partial}_2{\xi}^2)\]
From the mathematical point of view, one uses to consider the Lie operator ${\cal{D}}\xi={\cal{L}}(\xi)\omega:T\rightarrow S_2T^*$ (symmetric tensors), sometimes called Killing operator, through the formula:\\
\[      ({\cal{D}}\xi)_{ij}\equiv {\omega}_{rj}{\partial}_i{\xi}^r+{\omega}_{ir}{\partial}_j{\xi}^r+{\xi}^r{\partial}_r{\omega}_{ij}={\Omega}_{ij}=2{\epsilon}_{ij}   \]
One may check at once the only generating "{\it compatibility condition}" ${\cal{D}}_1\epsilon=0$, namely:\\
\[    {\partial}_{11}{\epsilon}_{22}+{\partial}_{22}{\epsilon}_{11}-2{\partial}_{12}{\epsilon}_{12}=0 \]
which is nothing else than the Riemann tensor of a metric, linearized at $\omega$.\\
However, the main experimental reason for introducing the first operator of this type of Janet sequence is the fact that the deformation is made from the displacement and first derivatives but must be invariant under any rigid motion. In the general case it must therefore have $(n+n^2)-(n+n(n-1)/2)=n(n+1)/2$ components, that is 3 when $n=2$, and this is the reason why introducing the deformation tensor $\epsilon$. For most finite element computations, the action density (local free energy) $w$ is a (in general quadratic) function of $\epsilon$ and people use to define the stress by the formula ${\sigma}^{ij}=\partial w/\partial {\epsilon}_{ij}$ which is not correct because $w$ only depends on ${\epsilon}_{11},{\epsilon}_{12}, {\epsilon}_{22}$ when $n=2$  as the deformation tensor is symmetric {\it by construction}. Finally, textbooks escape from this trouble by {\it deciding} that the stress should be symmetric and this is a {\it vicious circle} because we have proved it was not an assumption but an experimental result depending on specific constitutive laws. Accordingly, when $n=2$, we should have ${\sigma}^{ij}{\epsilon}_{ij} = {\sigma}^{11}{\epsilon}_{11}+(2{\sigma}^{12}){\epsilon}_{12}+{\sigma}^{22}{\epsilon}_{22}$. Hence, even if we find the correct stress equations with {\it this convenient duality} keeping the factor "2", we have no way to get the stress AND couple-stress equations TOGETHER.\\

\noindent
3){\it Only the Spencer sequence can be used}:\\
Let us construct the formal adjoint of the Spencer operator by multiplying all the $(2\times 2)+2=6$ linearly independent nonzero components by corresponding test functions. For simplifying the summation, we shall raise and lower the indices by means of the (constant) euclidean metric, setting in particular ${\xi}_i={\omega}_{ir}{\xi}^r$ and ${\xi}_{i,j}={\omega}_{ir}{\xi}^r_j$. Comparing to Example 8, the only nonzero first jets coming from the $2\times 2$ skewsymmetric infinitesimal rotation matrix of first jets are now ${\xi}_{1,2}=-{\xi}_{2,1}$ while the second order jets are zero because isometries are linear transformations. We obtain in the present situation:\\
\[ {\sigma}^{11}{\partial}_1{\xi}_1+{\sigma}^{12}({\partial}_1{\xi}_2-{\xi}_{1,2})+{\sigma}^{21}({\partial}_2{\xi}_1-{\xi}_{2,1})+{\sigma}^{22}{\partial}_2{\xi}_2+{\mu}^{r,12}{\partial}_r{\xi}_{1,2}  \]
Integrating by parts and changing the sign, we just need to look at the coefficients of ${\xi}_1,{\xi}_2$ and ${\xi}_{1,2}$, namely:\\
\[  \begin{array}{lcl}
 {\xi}_1 & \longrightarrow  &  {\partial}_1{\sigma}^{11}+{\partial}_2{\sigma}^{21}=f^1  \\
 {\xi}_2 & \longrightarrow   & {\partial}_1{\sigma}^{12}+{\partial}_2{\sigma}^{22}=f^2  \\
 {\xi}_{1,2} & \longrightarrow  & {\partial}_r{\mu}^{r,12}+{\sigma}^{12}-{\sigma}^{21}=m^{12}  
\end{array}  \]
in order to get the adjoint operator $ad(D):{\wedge}^{n-1}T^*\otimes R_1^*\rightarrow {\wedge}^nT^*\otimes R_1^*:(\sigma,\mu)\rightarrow (f,m)$ relating for the first time the torsor framework to the Lie coalgebroid $R_1^*$. These equations are {\it exactly} the three stress and couple-stress equations of 2-dimensional elasticity. In the n-dimensional case, a similar calculation, left to the reader as an exercise of indices, should produce {\it exactly} the $n(n+1)/2$ stress and couple-stress equations in general. It is now possible to enlarge the group in order to get more equations, that is {\it as many equations as the number of group parameters}. Using the conformal group of space-time, the 4 elations give rise to 4 nonzero second order jets only which allow to exhibit the 4  Maxwell equations for the induction $(\vec{H},\vec{D})$ along lines only sketched by H. Weyl in [20] because the needed mathematics were not available before 1970. But, as we already said, this is another story !.\\

\noindent
{\bf REMARK 7}: It becomes now clear that the $n^2(n^2-1)/4$ first order compatibility conditions for the Cosserat fields [11] (the so-called torsion and curvature of E. Cartan [4]) are described by the second Spencer operator in the Spencer sequence while the $n^2(n^2-1)/12$ second order compatibility conditions for the deformation tensor (the so-called Riemann curvature) are described by the second operator ${\cal{D}}_1$ in the Janet sequence. Accordingly, {\it the torsion+curvature of Cartan is not at all the generalization of the curvature of Riemann}, contrary to what is still claimed in mathematical physics today.\\

\noindent
{\bf CONCLUSION}:\\

The original approach of V. Arnold to hydrodynamics was based on specific analytic results only valid for the Lie pseudogroup of volume preserving transformations. We have extended this approach to an arbitrary Lie pseudogroup, meanwhile showing that the previous results are in fact purely formal results based on a new approach to duality theory in constrained variational calculus.\\
However, we have also found that these techniques, where the "time" variable is considered as a simple parameter, were not "natural" in some sense and could be superseded by the construction of the nonlinear Spencer sequence for an arbitrary Lie groupoid. This result, explaining the gauging concept in the jet framework, also achieves the modern interpretation of the works done at the beginning of the previous century by the brothers E. and F. Cosserat on the foundation of continuum mechanics and by H. Weyl on the foundation of electromagnetism.\\

\noindent
{\bf BIBLIOGRAPHY}:\\
\noindent
1) ARNOLD, V.: Sur la g\' eom\' etrie des groupes de Lie de dimension infinie et ses applications \`a l'hydrodynamique des fluides parfaits, Ann. Inst. Fourier (Grenoble), 16, 1, 1966, pp. 319-361.\\
\noindent
2) ARNOLD, V.: M\' ethodes math\' ematiques de la m\' ecanique classique, Appendice 2 (G\' eod\' esiques des m\' etriques invariantes \`a gauche sur des groupes de Lie et hydrodynamique des fluides parfaits), MIR, Moscow, 1974, 1976.\\
\noindent
3) BIRKHOFF, G.: Hydrodynamics, Princeton University Press, , Princeton, 1954; French translation: Hydrodynamique, Dunod, Paris, 1955.\\
\noindent
4) CARTAN, E.: Sur une g\' en\' eralisation de la notion de courbure de Riemann et les espaces \`a torsion, C. R. Acad\' emie des Sciences Paris, 174, 1922, P. 522.\\
\noindent
5) CHETAEV, N.G.,:C.R. Acad\' emie des Sciences, Paris, 185, 1927,p. 1577.\\
\noindent
6) CHETAEV, N.G.: Theoretical Mechanics, MIR, Moscow, 1989 and Springer Verlag.\\
\noindent
7) CHWOLSON,O.D.: Trait\' e de physique.(In particular III,2,p. 537 + III,3, p. 994 + V, p. 209), Hermann, Paris, 1914.\\
\noindent
8) COSSERAT, E. and F.: Note sur la th\' eorie de l'action euclidienne. In : Trait\' e de m\' ecanique rationelle (Appell, P., ed.),t. III, pp. 557-629, Gauthiers-Villars, Paris, 1909.\\
\noindent
9) COSSERAT, E. and F.: Th\' eorie des corps d\' eformables, Hermann, Paris, 1909.\\
\noindent
10) HERGLOTZ, G. : Uber die mechanik des deformierbaren k\" orpers vom standpunkt der relativit\" atstheorie, Ann. der Physik, 36, 1911,pp. 493-517.\\
\noindent
11) KOENIG, G.: Le\c{c}ons de cin\' ematique (The Note "Sur la cin\' ematique d'un milieu continu" by E. and F. Cosserat, pp. 391-417, has never been quoted elsewhere), Hermann, Paris, 1897.\\
\noindent
12) KUMPERA, A., SPENCER, D. C.: Lie equations, ANN. Math. Studies 73, Princeton University 
Press, Princeton, 1972.\\
\noindent
13) POINCARE, H.: Sur une forme nouvelle des \' equations de la m\' ecanique, C. R. Acad\' emie des Sciences Paris, 132, 7, 1901, p. 369-371. \\
\noindent
14) POMMARET, J. -F.: Lie pseudogroups and mechanics, Gordon and Breach, New York, 1988.\\
\noindent
15) POMMARET, J.-F.: Partial differential equations and Lie pseudogroups : New perspectives for applications, Kluwer, Dordrecht, 1994.\\
\noindent
16) POMMARET, J.-F.: Group interpretation of coupling phenomena, Acta Mechanica, 149, 2001,pp. 23-39.\\
\noindent
17) POMMARET, J.-F.: Partial differential control theory, Kluwer, Dordrecht, 2001.\\
\noindent
18) POMMARET, J.-F. : Fran\c{c}ois Cosserat et le secret de la th\' eorie math\' ematique de l'\' elasticit\' e, Annales des Ponts et Chauss\' ees, 82, 1997, pp. 59-66.\\
\noindent
19) TEODORESCU, P. P. : Dynamics of linear elastic bodies, Editura Academiei, Bucuresti, Romania; Abacus Press, Tunbridge, Wells, 1975.\\
\noindent
20) WEYL, H. : Space, time, matter, Springer, Berlin, 1918, 1958; Dover, 1952.\\
\noindent
21) YANG, C. N. : Magnetic monopoles, fiber bundles and gauge fields, Ann. New York Acad. Sciences, 294, 1977, p. 86.\\

\noindent
RESEARCH REPORT / RAPPORT DE RECHERCHE  CERMICS  n$^{\circ}$ 353 , received may  2007, modified for publication may 2008\\
(http://cermics.enpc.fr)\\

\end{document}